\documentclass{jnmp}
\usepackage{amsmath}
\JNMPnumberwithin{equation}{section}
\theoremstyle{plain}
\newtheorem{thm}{Theorem}[section]
\newtheorem{lem}{Lemma}[section]
\newtheorem{cor}{Corollary}[section]
\theoremstyle{definition}
\newtheorem{rem}{Remark}[section]
\newcommand{\B}{\mathbf}
\newcommand{\C}{\mathcal}
\newcommand{\D}{\mathbb}
\DeclareMathOperator{\dist}{dist}
\renewcommand{\Im}{\operatorname{Im}}
\DeclareMathOperator{\meas}{meas}
\renewcommand{\Re}{\operatorname{Re}}
\begin{document}
\renewcommand{\evenhead}{A Boutet de Monvel, E Khruslov and V Kotlyarov}
\renewcommand{\oddhead}{Soliton Asymptotics of Rear Part
                         of Non-Localized Solutions of the K-P Equation}
\thispagestyle{empty}
\FirstPageHead{9}{1}{2002}{\pageref{firstpage}--\pageref{lastpage}}{Article}
\copyrightnote{2002}{A Boutet de Monvel, E Khruslov and V Kotlyarov}
\Name{Soliton Asymptotics of Rear Part
       of Non-Localized Solutions of the Kadomtsev-Petviashvili Equation}
\label{firstpage}
\Author{Anne BOUTET de MONVEL~$^\dagger$, Eugene KHRUSLOV~$^{\ddagger\ast}$,
               and Vladimir KOTLYAROV~$^{\ddagger\star}$}
\Address{$^\dag$ Universit\'e Paris-7,
          Physique math\'ematique et G\'eom\'etrie,
          UFR de Math\'ematiques, \\
          ~~case 7012,
          2 place Jussieu, 75251 Paris Cedex 05, France\\
          ~~E-mail: aboutet@math.jussieu.fr\\[10pt]
          $^\ddag$ Mathematical Division,
          Institute for Low Temperature Physics,
          47 Lenin Avenue, \\
          ~~310164 Kharkov, Ukraine\\
          $^\ast$ E-mail: khruslov@ilt.kharkov.ua\\
          $^\star$ E-mail: vkotlyarov@ilt.kharkov.ua}
\Date{Received April 02, 2001; Revised November 03, 2001;
       Accepted November 04, 2001}
\begin{abstract}
\noindent
We construct non-localized, real global solutions of the 
Kadomtsev-Petviashvili-I
equation which vanish for $x\to-\infty$ and study their large time
asymptotic behavior. We prove that such solutions eject (for $t\to\infty$)
a train of curved asymptotic solitons which move behind the
basic wave packet.
\end{abstract}
\section{Introduction}

The study of long-time asymptotic behavior of solutions of nonlinear 
evolution  equations
has attracted growing attention in last years and a number of  papers 
have been devoted to
this problem (see review paper \cite{DIZ} and references therein). 
Interest in this
problem was especially stimulated by the discovery of the inverse 
scattering transform method
\cite{GGKM,ZMNP,AKNS,FT}. In particular, one remarkable result obtained by this
method was the proof that any localized (i.e.\ rapidly decreasing as 
$|x|\to\infty$)
solution $u(x,t)$ of the Korteweg-de Vries (KdV) equation
\[
u_t-6uu_x+u_{xxx}=0
\]
splits into a finite number of solitons when time tends to infinity \cite{Sha}
(see also \cite{T,Schuur}). Other nonlinear evolution equations with one space
variable integrable by the inverse scattering transform also exhibit similar 
splitting.
This phenomenon is an argument in favor of the physical 
interpretation of solitons as stable
``long-living'' particles.

The simplest non-localized solution is of step-like form, i.e.\ the 
solution with
following
asymptotic behavior:
\[
u(x,t)=
\begin{cases}
0,&x\to+\infty\\
-c^2,& x\to-\infty.
\end{cases}
\]
It was proved in \cite{Kh} that the step-like solution of the KdV 
equation splits
into $\bigl[\frac{N+1}{2}\bigr]$ soliton-like objects in the
neighbourhood
\[
G^+_N(t)=\{x\in\D{R}\mid x>4c^2t-N\ln t\}
\]
of solution front
($N\in\D{Z}^+$ and $[\;\cdot\;]$ denotes integer part).
The form of these objects
is similar to ordinary soliton but their velocities depend on $t$. In 
contrast with
ordinary solitons they are not exact solutions of the KdV equation, 
however they satisfy it
with increasing accuracy when
$t\to\infty$. For this reason such objects are called ``asymptotic 
solitons''. The number of
these asymptotic solitons increases to infinity when $t\to\infty$ if 
the observation domain
in the neighbourhood of the the solution front is extended correspondently.

The same phenomenon of generation of asymptotic solitons trains on 
the solution front
takes place (under certain conditions) for non-localized solutions of 
more general form as
well as for other nonlinear evolution equations integrable by the 
inverse scattering
transform \cite{KhK}. Roughly speaking this phenomenon can be 
considered as a manifestation
of the fact that any non-localized initial data consists of an 
infinite number of solitons
which are gradually ejected at the front (first is the most rapid of 
them). An important
condition for the ejection is the existence of wide ``living area'' 
for solitons where they
can propagate without collisions. In the case where the non-localized 
initial data vanish
for $x\to\infty$ this area is a positive beam $(at,\infty)$ for some 
suitable $a>0$.

Another type of non-localized solutions of the KdV equation which 
vanish as $x\to-\infty$
was considered in \cite{BarK}. It has been proved that under some 
conditions trains of
asymptotic solitons are formed on a tail of the solution as 
$t\to\infty$. These solitons
move to the right following behind the basic wave packet. Physically 
it can be treated as
an ejection of the slower solitons from a non-localized initial 
perturbation (initial data
$u(x,0)$). This phenomenon takes place under certain conditions on 
the spectrum of the
Schr\"odinger operator whose potential is the initial data $u(x,0)$.
This condition is the
existence of a gap between the continuous spectrum of multiplicity 
one, which is provided
by a non-trivial asymptotic behavior of $u(x,0)$ as $x\to\infty$, and 
the continuous
spectrum of multiplicity two (positive half-axis).

In \cite{AKhK1,A,AKhK2,OPKh} similar problems were studied for the
Kadomtsev-Petviashvili  equations (KP-I and KP-II):
\begin{equation} 
\label{KP}
\frac{\partial}{\partial x}\Bigl(u_t+\frac{3}{2}u u_x
+ \frac{1}{4}u_{xxx}\Bigr) +\frac{3}{4}\alpha^2u_{yy}=0,
\end{equation}
where $\alpha=i$ for KP-I and $\alpha=1$ for KP-II.
These are equations with two spatial variables $x, y$ integrable by 
inverse scattering
transform. See also \cite{IgABM-00,IgABM-01} for the Johnson equation 
and for the modified
Kadomtsev-Petviashvili-I equation (mKP-I).

Following the Zakharov-Shabat scheme of the dressing method we look for
solutions in the form:
\begin{equation} 
\label{u}
u(x,y,t)=2\frac{d}{d x}K^\pm(x,x,y,t),
\end{equation}
where the function $K^\pm(z,x,y,t)$ is a solution of the Marchenko 
integral equation:
\begin{equation} 
\label{Meq}
K^\pm(z,x,y,t)+F(z,x,y,t)\pm
\int_x^{\pm\infty} K^\pm(s,x,y,t)F(z,s,y,t)\,d s=0
\end{equation}
viewed as an equation in $z>x$ ($z<x$) with parameters $x,y,t$.
The kernel $F(z,x,y,t)$ in (\ref{Meq}) satisfies the system of linear 
differential
equations:
\begin{equation} 
\label{Feq}
\begin{cases}
F_t+F_{xxx}+F_{zzz}=0,&\\
\alpha F_y+F_{xx}-F_{zz}=0.&
\end{cases}
\end{equation}
Choosing solutions of this system in an appropriate way and solving 
the integral equation
(\ref{Meq}) (with the sign ``$+$'' or ``$-$'') one can construct 
solutions of the KP
equations which vanish as $x\to+\infty$ or $x\to-\infty$. For the 
KP-I equation global
real solutions vanishing as $x\to+\infty$ were constructed in 
\cite{A} (see also
\cite{AKhK2}). It was proved that in the neighbourhoods of their fronts:
\[
x>C(Y)t-\frac{1}{a(Y)}\log t^N,
\]
these solutions are asymptotically represented as
follows:
\begin{align}	                                                \label{usol}
u(x,y,t)&=\sum_{n=1}^N
2a^2(Y)\left(\cosh\Bigl[a(Y)(x-C(Y)t+\frac{1}{a(Y)}\log t^{n-1/2})-
\log\varphi_n(Y)\Bigr]\right)^{-1}\notag\\
&\quad
+\mathrm{O}(t^{-1/2+\varepsilon})
\end{align}
where $N$ is any natural number, $Y=y/t$, $a(Y)$, $C(Y)$ and 
$\varphi_n(Y)$ are some
positive functions. This means that the trains of curved solitons are 
formed in the
neighbourhood of the solution front because the ridges of these 
solitons at time $t$ are
located along curves $x=C(y/t)t-a^{-1}(y/t)\log t^{n-1/2}+ \log
\varphi_n(y/t)$ ($n=1,2,\dots$). The phenomenon of formation of curved solitons
trains nearly solution fronts is observed also for the KP-II equation, however
in this case even solitons have a singularity \cite{AKhK1}.

This paper is devoted to the study of non-localized solutions of the 
KP-I equation which
vanish as $x\to-\infty$. We construct a global real solution of this 
type and prove that a
back part of the solution splits into curved asymptotic solitons. The 
dressing method of
Zakharov-Shabat described above is not suitable for the investigation 
of long-time
asymptotics of the tail of solutions. Therefore we use another method 
based on an integral
equation in the plane of spectral parameters. This approach is also 
suitable for the study
of long-time asymptotics of non-localized solutions in the 
neighbourhood of the solution
front, which was first shown in \cite{KhS} for the KdV equation.

The paper is organized as follows. In Section 2 we prove the 
existence of a global real
solution of the KP-I equation, vanishing as $x\to-\infty$. In Section 
3 we reduce the
problem to a degenerated integral equation and obtain a determinant 
formula for the solution.
In Section 4 we study the asymptotic behavior of the determinant 
formula as $t\to\infty$
and prove an ejection of curved solitons which move behind the basic 
wave packet.

\section{Construction of global solutions vanishing as $x\to-\infty$}

Let $\Omega$ be a domain in the complex plane $\D{C}$ ($k=p+iq$) with 
smooth boundary
$\Gamma=\partial\Omega$ located in the right half plane at positive 
distance from the
imaginary axis. Let us define the function $E(p,q):=E(p,q,x,y,t)$
on $\Omega$ as follows:
\begin{equation} 
\label{E}
E(p,q)=e^{p(x-f(p,q,Y)t)},
\end{equation}
with
\begin{equation} 
	\label{f}
f(p,q,Y)= p^2-3q^2-2qY,
\end{equation}
where $x,y,t\in\D{R}$, and $Y=y/t$ are parameters. Let us now 
consider the integral
equation (with respect to $\psi(p,q):=\psi(p,q,x,y,t)$):
\begin{equation} 
\label{psieq}
\psi(p,q)+\int_\Omega\frac1{\lambda+\overline{k}}\;{E(p,q)E(\mu,\nu)}
\psi(\mu,\nu)g(\mu,\nu)\,d\mu\,d\nu=E(p,q),
\end{equation}
where $\lambda=\mu+i\nu\in\Omega$,\ $k=p+iq\in\Omega$, $\bar k=p-iq$ 
and $g(\mu,\nu)$ is a
smooth positive function on $\overline\Omega$. Furthermore, if 
$\Omega$ is unbounded,
$g(\mu,\nu)$ satisfies
\begin{equation} 
	\label{g}
\int_\Omega e^{c\mu(\mu^2+\nu^2)}g(\mu,\nu)\,d\mu\,d\nu<\infty
\end{equation}
for any $c>0$.

\begin{lem}                                                       \label{lem.1}
There exists a unique solution $\psi(p,q)=\psi(p,q,x,y,t)$
of $(\ref{psieq})$ which is $C^{\infty}$ with respect to
$x,y,t\in\D{R}$.
\end{lem}

\begin{proof}
${L}^2_g(\Omega)$ denotes the Hilbert space of complex valued functions
on $\Omega$ with norm
\[
||\varphi||:=\left\{\int_\Omega|\varphi(\mu,\nu)|^2
g(\mu,\nu)d\mu\,d\nu\right\}^{\frac{1}{2}}.
\]
Let $\B{A}$ be an operator on $L^2_g(\Omega)$, depending on 
parameters $x,y,t\in\D{R}$,
as follows:
\begin{equation} 
\label{A}
[\B{A}\varphi](p,q)=\int_\Omega\frac{E(p,q)E(\mu,\nu)}{\lambda+\bar k}
\varphi(\mu,\nu)g(\mu,\nu)\,d\mu\,d\nu.
\end{equation}
According to (\ref{E}), (\ref{f}) and (\ref{g}), $\B{A}$ is Hilbert-Schmidt
for any values of the parameters $x,y,t$. For $\varphi\in L^2_g(\Omega)$ we
can write
\[
(\B{A}\varphi,\varphi)=\int_\Omega E(p,q)\int_\Omega
\frac{E(\mu,\nu)}{\lambda+\bar k}\varphi(\mu,\nu)g(\mu,\nu)\,d\mu\,d\nu\,
\bar\varphi(p,q)g(p,q)\,d p\,d q,
\]
where $(\,\cdot\,,\,\cdot\,)$ denotes the inner product in $L^2_g(\Omega)$ and
$\lambda=\mu+i\nu$, $\bar k=p-i q$. Using the equality
\[
\frac{1}{\lambda+\bar k}=\int^0_{-\infty}e^{(\lambda+\bar k)s}\,d s,
\]
which is true because $\Re\lambda>0$, $\Re k>0$, we obtain
\begin{equation} 
\label{A>0}
(\B{A}\varphi,\varphi)=
\int^0_{-\infty}d s
\left|\int_\Omega
E(\mu,\nu)\varphi(\mu,\nu)e^{(\mu+i\nu)s} g(\mu,\nu)\,d\mu\,d\nu\right|^2
\geq 0.
\end{equation}
It follows from (\ref{A>0}) that the operator $\B{A}$ is positive. Therefore,
the homogeneous equation $\varphi+\B{A}\varphi=0$ has only the trivial solution
$\varphi\equiv0$. Since $\B{A}$ is Hilbert-Schmidt, hence compact in
$L^2_g(\Omega)$, the inhomogeneous equation
\begin{equation} 
\label{psi'}
\psi+\B{A}\psi=\zeta
\end{equation}
has a unique solution in $L^2_g(\Omega)$ for any
$\zeta\in L^2_g(\Omega)$. Let us note now that the integral equation 
(\ref{psieq}) has the
same form as (\ref{psi'}) with $\zeta=E(p,q)=E(p,q,x,y,t)$ belonging 
to the space
$L^2_g(\Omega)$ due to (\ref{E}). Therefore, (\ref{psieq}) has a 
unique solution
$\psi(p,q)$, depending on the parameters $x,y,t$. The first 
derivatives $D\psi$ of this
solution with respect to $x,y,t$ also satisfy (\ref{psi'}) with right-hand side
\[
\zeta=DE(p,q)-\int_\Omega D\,\frac{E(p,q)E(\mu,\nu)}{\lambda+\bar k}
\psi(\mu,\nu)g(\mu,\nu)\,d\mu\,d\nu
\]
belonging to the space $L^2_g(\Omega)$. This proves their existence.
Existence of high order derivatives is proved by induction.
\end{proof}

\begin{cor}
The inverse operator $(\B{I}+\B{A})^{-1}$ exists and its norm in
$\C{L}(L^2_g(\Omega))$
is uniformly bounded with respect to $x,y,t\in\D{R}$:
\begin{equation} 
	\label{A<1}
||(\B{I}+\B{A})^{-1}||\le 1.
\end{equation}
\end{cor}

Let us now define the function
\begin{equation}	\label{u2}
u(x,y,t)=-2\frac{\partial}{\partial x}\int_\Omega
E(p,q,x,y,t)\psi(p,q,x,y,t)g(p,q)\,d p\,d q,
\end{equation}
where $\psi(p,q,x,y,t)$ is the solution of the integral equation 
(\ref{psieq}), and
$E(p,q,x,y,t)$ is determined by $(\ref{E})$ and $(\ref{f})$.

\begin{lem} 
\label{lem.2}
The function $(\ref{u2})$ is a solution of the KP-I equation
$(\ref{KP})$.
It is a real solution, vanishing as $x\to-\infty$.
\end{lem}

\begin{proof}
Let us multiply (\ref{psieq}) by $E(p,q,z,y,t)e^{i q(z-x)}g(p,q)$ and
integrate with respect to $p,q\in\Omega$. Then, using
\[
\frac{1}{\lambda+\bar k}=\int^x_{-\infty}e^{(\lambda+\bar k)(s-x)}\,d s,
\]
together with (\ref{E}), (\ref{f}) we obtain
\begin{equation} 
	\label{M}
K(z,x,y,t)+\int_{-\infty}^xF(z,s,y,t)K(s,x,y,t)d s+F(z,x,y,t)=0,
\end{equation}
where
\[
K(z,x,y,t)=-\int_\Omega
E(p,q,z,y,t)e^{i q(z-x)}\psi(p,q,x,y,t)g(p,q)\,d p\,d q,
\]
\[
F(z,x,y,t)=\int_\Omega
E(p,q,z,y,t)E(p,q,x,y,t)e^{i q(z-x)}g(p,q)\,d p\,d q.
\]
Taking into account (\ref{E}), (\ref{f}), and (\ref{g}) it is easy to
show that the function $F(z,x,y,t)$ satisfies equations (\ref{Feq}). 
It is also obvious
that formula (\ref{u2}) and equation (\ref{M}) (with respect to
$K(z,x,y,t)$,
$z<x$) correspond to the equations (\ref{u}), (\ref{Meq}) (with the 
``$-$'' sign). Thus
according to the Zakharov-Shabat dressing method, $u$ defined by 
(\ref{u2}) is a
solution of the KP-I equation. Taking into account that the domain 
$\Omega$ is contained in
the right half plane at positive distance from the imaginary axis and 
that the function
$g(p,q)$ satisfies inequality (\ref{g}) it is easy to prove that 
solution (\ref{u2})
vanishes for $x\to-\infty$.

This solution is real. Indeed according to (\ref{u2})
it is sufficient to prove $\Im(E,\bar\psi)=0$, where $\psi$ is solution of
the integral equation
(\ref{psieq}) and $(\,\cdot\,,\,\cdot\,)$ denotes the inner product 
in $L^2_g(\Omega)$.
Let us
remind that equation (\ref{psieq}) in $L^2_g(\Omega)$ takes the form 
(\ref{psi'})
with $\B{A}$ a positive operator and real right hand side $\zeta=E$. 
Then we obtain
\[
(E,\bar\psi)=\overline{(\zeta,\psi)}=\overline{(\psi+\B{A}\psi,\psi)}=
(\psi,\psi)+(\B{A}\psi,\psi)\geq 0.
\]
Hence $\Im(E,\bar\psi)=0$ and the lemma is proved.
\end{proof}

Thus for a given function $g(p,q)\ge0$, formulas (\ref{u2}) and
(\ref{E})-(\ref{psieq}) represent a global real solution $u(x,y,t)$ of the KP-I
equation. This solution decays as $x\to-\infty$, but its behavior as 
$x\to+\infty$ is
unknown.

Our main goal is to study the asymptotic behavior of this solution as 
$t\to\infty$
in suitable neighbourhoods of the rear part of the solution, namely 
neighbourhoods of the
form:
\[
D_N(t)=\Bigl\{(x,y)\,\Bigm|\,-\infty<y<\infty,\ x<C(Y)t+ 
\frac{1}{a(Y)}\log t^N\Bigr\}.
\]
We show that in such
domains and for large $t$ ($t>T(N)$) the solution $u(x,y,t)$ has the 
asymptotic behavior
described by (\ref{usol}), where
$a(Y)$, $C(Y)$, $\varphi_n(Y)$ are expressed in terms of $g(p,q)$.
This formula means that
the solution splits into a sequence of curved solitons, which are formed
in the neighbourhood
of its trailing edge. To prove this asymptotic formula we first need
to approximate the
solution of the integral equation (\ref{psieq}) by solutions of 
integral equations with
appropriate degenerate kernels.

\section{Integral equation with degenerate kernel}

Let $k_0=p_0+i q_0\in\Gamma$ be an arbitrary point of the boundary of $\Omega$.
We will use the following double power series expansion:
\begin{equation}                                                  \label{cij}
\frac{1}{\lambda+\bar k}=
\sum_{i,j=0}^\infty C_{ij}(\lambda-k_0)^i (\bar k-\bar k_0)^j,
\end{equation}
where
\[
C_{ij}=(-1)^{i+j}\frac{(i+j)!}{i!j!(2p_0)^{i+j+1}}.
\]
It is easy to check that (\ref{cij}) converges in the polydisk
\[
\Pi(k_0)=\{(\lambda,k)\,\mid\, |\lambda-k_0|<p_0,\, |k-k_0|<p_0\}.
\]
Below we will choose $k_0$ as a point of $\overline\Omega$ where the
function $f(p,q,Y)=p^2-3q^2-2qY$ attains its minimal value.
We suppose such a point exists
and is unique. This point depends on parameter $Y$:
\begin{equation} 
	\label{k0}
k_0=k_0(Y)=p_0(Y)+iq_0(Y).
\end{equation}
Let us denote by $C(Y)$ the value of the function $f(p,q,Y)$
at the point $k_0(Y)$, i.e.
\begin{equation}                                                 \label{CY}
C(Y)=f(p_0(Y),q_0(Y),Y)=\min_{(p,q)\in\Omega}f(p,q,Y)
\end{equation}
and by $\chi_{N,Y}(p,q)$ the characteristic function of a subdomain
$G_{N,Y}\subset\Omega$ such that
\begin{equation}                                                 \label{GN}
0<\dist(k_0(Y),
\Omega\setminus\overline G_{N,Y})<\frac{p_0(Y)}{2}.
\end{equation}
This subdomain depends on $Y$ and on $N$. It will be precisely defined later.

Using the expansion (\ref{cij}) we represent the operator
$\B{A}$ (\ref{A}) in the form:
\begin{equation}                                                 \label{ABC}
\B{A}=\B{A}_N+\B{B}_N+\B{C}^1_N+\B{C}^2_N,
\end{equation}
where the operators $\B{A}_N$, $\B{B}_N$, $\B{C}^1_N$, and 
$\B{C}^2_N$ are defined by
\begin{align*}
[\B{A}_N\varphi](p,q)
&=
\int_{\Omega}E(p,q)E(\mu,\nu)\chi_{N,Y}(p,q)
\chi_{N,Y}(\mu,\nu)g(\mu,\nu)\\
&\qquad
\times\sum_{i,j=0}^N C_{ij}(\lambda-k_0)^i(\overline
k-\overline k_0)^j\varphi(\mu,\nu)\,d\mu\,d\nu,\\
[\B{B}_N\varphi](p,q)
&=
\int_{\Omega}E(p,q)E(\mu,\nu)\chi_{N,Y}(p,q)\chi_{N,Y}(\mu,\nu)g(\mu,\nu)\\
&\qquad
\times\sum_{(i,j)\in\tilde R^{(N)}} C_{ij}(\lambda-k_0)^i(\overline k-\overline
k_0)^j\varphi(\mu,\nu)\,d\mu\,d\nu,\\
[\B{C}^1_N\varphi](p,q)
&=
\int_{\Omega} \frac{E(p,q)E(\mu,\nu)}{\lambda+ \overline k} 
(1-\chi_{N,Y}(\mu,\nu))
g(\mu,\nu)\varphi(\mu,\nu)\,d\mu\,d\nu,\\
[\B{C}^2_N\varphi](p,q)
&=
\int_{\Omega} \frac{E(p,q)E(\mu,\nu)}{\lambda+ \overline k} (1-\chi_{N,Y}(p,q))
\chi_{N,Y}(\mu,\nu)g(\mu,\nu)\varphi(\mu,\nu)\,d\mu\,d\nu.
\end{align*}
Here $\lambda= \mu+i\nu$, $k=p+iq$ and
\begin{align*}
\tilde R^{(N)}
&:=\{(i,j)\mid i,j\geq 0\}
\setminus\{(i,j)\mid 0\le i,j\le N\}\\
&\phantom{:}=
\{(i,j)\mid i\geq 0,\, j\geq N+1\}\cup
\{(i,j)\mid i\geq N+1,\, j\geq 0\}.
\end{align*}
Let us estimate the norm of the operators $\B{B}_N$, $\B{C}^1_N$ and 
$\B{C}^2_N$
in the space $L^2_g(\Omega)$. To avoid unessential complications we assume that
$\Omega$
is bounded. We will denote
\begin{align*}
&a=\inf_{\lambda\in\Omega}\Re\lambda,\qquad\qquad\qquad\qquad
b=\sup_{\lambda\in\Omega}\Re\lambda,\\
&d(\xi)= (b+a)\xi+(b-a)|\xi|,\qquad
\xi=x-C(Y)t,
\end{align*}
where $C(Y)$ is defined by (\ref{CY}). The norms of the operators
$\B{C}^i_N$ ($i=1,2$) on $L^2_g(\Omega)$ can be estimated as follows:
\begin{align*}
||\B{C}^i_N||^2
&\leq
\iint_{\Omega\times\Omega}(|\lambda+\overline k|)^{-2}
e^{2(p+\mu)\xi}\,e^{2p(C(Y)-f(p,q,Y))t}\,e^{2\mu(C(Y)-f(\mu,\nu,Y))t}\\
&\qquad\qquad
\times(1-\chi_{N,Y}(\mu,\nu))\;g(p,q)g(\mu,\nu)\,d p\,d q\,d\mu\,d\nu\\
&\le
\frac{e^{2d(\xi)-2am_0t}}{(2a)^2}\,\hat g^2(\meas\Omega)^2,
\end{align*}
where
\[
\hat g=\max_{(p,q)\in\Omega}g(p,q) \text{ and }
m_0=\min_{(p,q)\in\Omega\setminus G_{N,Y}}[f(p,q,Y)-C(Y)].
\]
According to (\ref{f}), (\ref{CY}), (\ref{GN}),
$m_0>0$. Hence,
\begin{equation}                                                  \label{CiN}
||\B{C}^i_N||\leq\hat g\times\meas\Omega\times e^{-am_0t/2}
\end{equation}
for $\xi<\frac{am_0t}{4b}$.
Taking into account that
$\tilde R^{(N)}\subset\cup_{k\geq N+1}\{(i,j)\mid i+j=k,\, i,j\ge0\}$
we can write (for $\lambda, k\in G_{N,Y}$):
\begin{align*}
\sum_{(i,j)\in\tilde R^{(N)}}|C_{ij}|\cdot|\lambda-k_0|^i |\overline 
k-\overline k_0|^j
&\le
\sum_{\substack{i,j\ge0,\\
                 i+j\ge N+1}}
|C_{ij}|\cdot|\lambda-k_0|^i |k-k_0|^j\\
& 
=\sum^\infty_{l=N+1}\sum_{i=0}^l\frac{l!}{i!(l-i)!(2p_0)^{l+1}}|\lambda-k_0|^i\, 
|k-k_0|^{l-i}\\
&=
\frac{1}{2p_0}\sum_{l=N+1}^\infty
\left(\frac{|\lambda-k_0|+|k-k_0|}{2p_0}\right)^l\\
&\leq\frac{1}{p_0}\left(\frac{|\lambda-k_0|+|k-k_0|}{2p_0}\right)^{N+1}.
\end{align*}
Using this inequality we obtain
\begin{align} 
\label{eBN}
||\B{B}_N||^2
&\leq\iint_{\Omega\times\Omega}
e^\Phi\,\frac{1}{p_0^2}\left(\frac{|\lambda-k_0|+|k-k_0|}{2p_0}\right)^{2(N+1)}
\chi\,d p\,d q\,d\mu\,d\nu\notag\\
&=\frac{1}{p_0^2}\sum_{\substack{i,j\ge0,\\
                 i+j=2(N+1)}}
                 |C_{ij}|J_i(x,Y,t)J_j(x,Y,t),
\end{align}
where
\begin{align}                                                        \label{Ji}
& \Phi=2(p+\mu)[x-t(f(p,q,Y)-f(\mu,\nu,Y))] \notag\\
& \chi=\chi_{N,Y}(p,q)\chi_{N,Y}(\mu,\nu) g(p,q)g(\mu,\nu)\notag\\
& J_i(x,Y,t)=\int_{G_{N,Y}}e^{2p(x-f(p,q,Y)t)}|k-k_0|^i\,g(p,q)\,d p\,d q,
\end{align}
and the numbers $C_{ij}$ are those introduced in (\ref{cij}).

We now define more precisely the subdomain $G_{N,Y}\subset\Omega$.
We will suppose that $\Gamma$ is defined by
\[
\Gamma=\{(p,q)\mid\varphi(p,q)=0\}
\]
where $\varphi(k)=\varphi(p,q)\in C^2(\bar{\Omega})$ and that the
curvature of $\Gamma$ is everywhere positive.
Since the hyperbola
\[
H(Y):=\{(p,q)\mid f(p,q,Y)-C(Y)=0\}
\]
is tangent to $\Gamma$ at the point
$k_0=k_0(Y)\in\partial\Omega$ we can introduce in some neighbourhood 
of $k_0$ new
coordinates
\begin{align*}
r&:=2p(f(p,q,Y)-C(Y))=: F(p,q)\\
s&:=\Bigl(\frac{\partial\varphi}{\partial q}(k_0)(p-p_0)-
\frac{\partial\varphi}{\partial p}(k_0)(q-q_0)\Bigr) ||
\nabla\varphi(k_0)||^{-1}=: \Phi(p,q).
\end{align*}
It is evident that $s$ is the projection of $(p-p_0,q-q_0)$ on the tangent
to the boundary $\Gamma$ at $k_0$. It is easy to check that the 
equation of $\Gamma$ near
$k_0\in\Gamma$ is of the form
\[
r=\alpha_0s^2+\mathrm{O}(s^3),
\]
where
\begin{equation} 
\label{curvat}
\alpha_0=||\nabla F(k_0)||\,\frac{\hat\kappa_0\mp\kappa_0}{2}.
\end{equation}
$\hat\kappa_0$ and $\kappa_0$ are the curvatures of $\partial\Omega$ and
$H(Y)$ at the point $k_0$ respectively. The minus sign occurs when $\Gamma$ and
$H(Y)$ are on the same side of their common tangent. The plus sign 
occurs otherwise.

In any case $\alpha_0>0$. Hence, $\Gamma$ is given by
\begin{align*}
&s=\hat s_{\pm}(r),\qquad\quad r\geq 0,\\
&\hat s_{\pm}(r)=\pm\sqrt{r/\alpha_0}+\mathrm{O}(r).
\end{align*}
\begin{figure}[ht]
\unitlength 1mm
\begin{picture}(150,75)
\linethickness{1.4pt}
\bezier{260}(99,50)(125,29.33)(99,11)
\bezier{108}(99,50)(86.67,57.33)(76,51)
\bezier{108}(99,11)(86.67,3.33)(76,10)
\bezier{76}(76,51)(69,48)(66,37)
\bezier{72}(66,37)(63.67,29.33)(67,20)
\bezier{56}(67,20)(70.67,13)(76,10)
\linethickness{0.4pt}
\put(7,30){\vector(1,0){130}}
\put(20,0){\vector(0,1){77}}
\bezier{212}(60,33)(80,48.33)(108,50)
\bezier{204}(59,41)(76.67,56.33)(104,58)
\put(79,44){\line(0,1){0}}
\put(79,44){\line(0,0){0}}
\put(79,45){\line(0,0){0}}
\put(79,45){\line(0,0){0}}
\put(74,49){\line(0,0){0}}
\put(20,53){\line(1,0){60}}
\put(79.8,30){\line(0,1){22.3}}
\bezier{60}(66.8,39)(67.2,38.6)(67.6,38.2)
\bezier{60}(67.6,41)(68.45,40.25)(69.3,39.3)
\bezier{100}(68.7,43.3)(70.15,41.85)(71.6,40.4)
\put(70,45.3){\line(1,-1){3.76}}
\put(71.1,46.9){\line(1,-1){4.55}}
\put(72.7,48.4){\line(1,-1){5.05}}
\put(74,49.8){\line(1,-1){5.71}}
\put(76,51){\line(1,-1){6.05}}
\put(78,52){\line(1,-1){6.27}}
\put(80,52.7){\line(1,-1){6.32}}
\put(82,53.5){\line(1,-1){6.48}}
\put(84,54){\line(1,-1){6.5}}
\put(86,54){\line(1,-1){6.1}}
\put(88,54){\line(1,-1){5.75}}
\put(90.8,53.4){\line(1,-1){4.85}}
\put(93.4,53){\line(1,-1){4.1}}
\bezier{60}(96.8,51.1)(97.8,50.1)(98.8,49.1)
\put(95,63.2){\vector(0,-1){6.2}}
\put(91.05,40.4){\vector(-2,3){5.5}}
\put(120,40){\vector(-3,-1){9}}
\put(18.5,73){\makebox(0,0)[cc]{$q$}}
\put(133,32){\makebox(0,0)[cc]{$p$}}
\put(87,14){\makebox(0,0)[cc]{$\Omega$}}
\put(121.8,41.3){\makebox(0,0)[cc]{$\Gamma$}}
\put(94.6,38){\makebox(0,0)[cc]{$G_{N,Y}$}}
\put(79,52){\small$\bullet$}
\put(94,65){\makebox(0,0)[cc]{$H(Y)$}}
\put(109.5,59){\makebox(0,0)[cc]{$r=0$}}
\put(117.5,50.7){\makebox(0,0)[cc]{$r=\varepsilon_0(N)$}}
\put(78.7,56){\makebox(0,0)[cc]{$k_0(Y)$}}
\put(78.7,27.6){\makebox(0,0)[cc]{$p_0(Y)$}}
\put(15,53.5){\makebox(0,0)[cc]{$q_0(Y)$}}
\end{picture}
\caption{$\Omega$, $\Gamma$ and the subdomain $G_{N,Y}$}
\label{fig-1}
\end{figure}
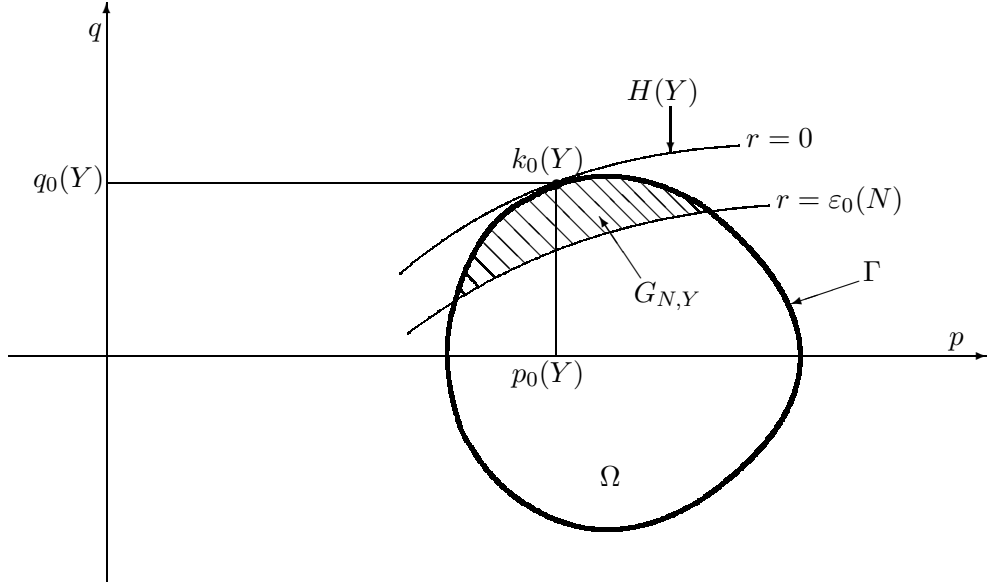
Now, the subdomain $G_{N,Y}\subset\Omega$ is defined as follows 
(Fig.~\ref{fig-1}):
\[
G_{N,Y}=\{(r,s)\mid 0<r<\varepsilon_0(N),\,\hat s_-(r)<s<\hat s_+(r)\}
\]
where $\varepsilon_0(N)$ is a small positive number, such that
\begin{equation} 
\label{epsilon}
\varepsilon_0(N)\le\frac{\alpha_0p_0^2}{16N^2
\bigl(\bigl|\frac{\partial p}{\partial
s}(k_0)\bigr|^2+\bigl|\frac{\partial p}{\partial r}(k_0)\bigr|^2\bigr)}.
\end{equation}
For $\varepsilon_0(N)>0$
small enough both bounds in (\ref{GN}) are clearly fulfilled.

Let us pass to an estimation of the integrals $J_i(x,Y,t)$ ($i=1,2$). 
Since for $k\in
G_{N,Y}$
\[
|k-k_0|=\sqrt{s^2+\frac{r^2}{||\nabla F(k_0)||^2}} 
\left(1+\mathrm{O}(r+|s|)\right)
\]
we have
\[
|k-k_0|^i\le 2^i\left(|s|^i+\frac{r^i}{|| \nabla F(k_0)||^i}\right).
\]
Due to (\ref{Ji}) we obtain
\begin{align}                                                 \label{Jib}
J_i(x,Y,t)
&\le 2^ie^{2p_0\xi}\int_0^{\varepsilon_0(N)}e^{-rt}\,d r\notag\\
&\quad
\times\int_{\hat s_-(r)}^{\hat s_+(r)}e^{2(p-p_0)\xi}\left(|s|^i+\frac{r^i}
{||\nabla F(k_0)||^i}\right)\hat g(r,s)|w(r,s)|\,d s,
\end{align}
where
$\xi=x-C(Y)t$, $\hat g(r,s)=g(p(r,s),q(r,s))$, and
$w(r,s)=\ds\frac{\partial(p,q)}{\partial(r,s)}$ is the Wronskian.
Using Taylor's series one can write:
\begin{equation}                                                 \label{exp}
e^{2(p-p_0)\xi}=1+2\left(\frac{\partial
p}{\partial s}(k_0)+ \frac{\partial p}{\partial r}(k_0)\right)E_0(r,s,\xi)
\end{equation}
where $E_0(r,s,\xi)$ is bounded:
\begin{equation}                                                 \label{E0}
|E_0(r,s,\xi)|\le
|\xi|e^{2\delta_0(N)|\xi|}
\end{equation}
for $(r,s)\in G_{N,Y}$, small
$\varepsilon_0(N)$ and
\begin{equation}                                                 \label{delta}
\delta_0(N):=2||\nabla p(k_0)||
\sqrt{\frac{\varepsilon_0(N)}{\alpha_0}}.
\end{equation}
Using this equality and the relations
\[
g(k)=g(k_0)+\mathrm{O}(r+|s|)\qquad w(k)=w(k_0)+\mathrm{O}(r+|s|)
\]
we obtain after integration of (\ref{Jib}) over $s$:
\begin{align}                                                     \label{Jib1}
J_i(x,Y,t)
&\le
\left(\frac{4}{\alpha_0}\right)^{\frac{i+1}{2}}\frac{g_0|w_0|}{i+1}e^{2p_0\xi}
\int_0^{\varepsilon_0(N)}
e^{-rt}r^{\frac{i+1}{2}}(1+\mathrm{O}(\sqrt r))\,d r\notag\\
&\quad
+|\xi|K_i(Y)e^{2p_0\xi+2\delta_0(N)|\xi|}\int_0^{\varepsilon_0(N)}
e^{-rt}r^{\frac{i+2}{2}}(1+\mathrm{O}(\sqrt r))\,d r\notag \\
&\le
K_{1i}(Y)\frac{e^{2p_0\xi}}{t^{\frac{i+3}{2}}}+
K_{2i}(Y)\frac{e^{2p_0\xi+2\delta_0(N)|\xi|}(1+|\xi|)}{t^{\frac{i+4}{2}}}
\end{align}
where $K_{1i}$ and $K_{2i}$ do not depend on $\xi$, $t$, and $g_0=g(k_0)$,
$w_0=w(k_0)$. Here
we have used the estimate
\begin{equation} 
	\label{gam}
\int_0^{\varepsilon_0(N)}e^{-rt}r^ddr<\int_0^\infty e^{-rt}r^ddr=
\frac{\Gamma(d+1)}{t^{d+1}},
\end{equation}
where $\Gamma(d)$ is the Euler $\Gamma$-function.
 From (\ref{eBN}) and (\ref{Jib1}) we derive the following bound for
the norm of
$\B{B}_N$:
\[
||\B{B}_N||\le K_N(Y)\left[\frac{e^{2p_0\xi}}
{t^{\frac{N+4}{2}}}+\frac{e^{2p_0\xi+\delta_0(N)|\xi|}(1+|\xi|)}
{t^{\frac{N+4}{2}+\frac{1}{4}}}+\frac{e^{2p_0\xi+2\delta_0(N)|\xi|}
(1+|\xi|)^2}{t^{\frac{N+5}{2}}}\right].
\]
Therefore, in view of (\ref{epsilon}) and (\ref{delta}), we find
\begin{equation}                                                 \label{bBN}
||\B{B}_N||\le\frac{K_N(Y)}{t^{1/2-\varepsilon}}
\qquad (0<\varepsilon<1/4).
\end{equation}

Now let us return to equation (\ref{psieq}) (or (\ref{psi'})).
Taking into account (\ref{A}) and (\ref{ABC}) we rewrite
it in the form:
\begin{equation}                                                 \label{psieq3}
\psi=\B{A}_N\psi+\B{B}_N^\prime\psi=e
\end{equation}
where $\B{B}_N^\prime=\B{B}_N+\B{C}^1_N+\B{C}^2_N$ and $e=E(p,q)$.
According to (\ref{CiN}) and (\ref{bBN}) the norm of
$\B{B}_N^\prime$ has the estimate
\begin{equation}                                                 \label{bBN'}
||\B{B}_N^\prime||\le
\frac{K_N^\prime(Y)}{t^{1/2-\varepsilon}}\quad
\text {if} \quad\xi=x-C(Y)t<\frac{1}{2p_0}\log t^{\frac{N+3}{2}+\varepsilon}.
\end{equation}
Let us look for the solution of equation (\ref{psieq3})
in the form
\begin{equation} 
\label{psi+delta}
\psi =\psi_N+\delta_N,
\end{equation}
where $\psi_N$ is the solution of the equation
\begin{equation}                                                 \label{psiN}
\psi_N+\B{A}_N\psi_N=e,
\end{equation}
and therefore
\begin{equation}                                                 \label{deln}
\delta_N=-(\B{I}+\B{A})^{-1}\B{B}_N^\prime\psi_N.
\end{equation}
According to (\ref{u}) and (\ref{psi+delta}) the solution of the KP-I equation
is represented in the form:
\begin{equation}                                                 \label{u3}
u(x,t)=-2\frac{d}{d x}(\psi,e)=-2\frac{d}{d x}(\psi_N,e)
- 2\frac{d}{d x}(\delta_N,e).
\end{equation}
Taking into account (\ref{deln}), the fact that $\B{A}$ is self-adjoint,
and relations (\ref{psieq3}), (\ref{psiN}) we can write
\begin{align}                                                    \label{dele}
(\delta_N,e)&=\notag
\left((\B{I}+\B{A})^{-1}\B{B}_N^\prime\psi_N,e\right)=
\left(\B{B}_N^\prime\psi_N,(\B{I}+\B{A})^{-1}e\right)\\ \notag
&=(\B{B}_N^\prime\psi_N,\psi)=(\B{B}_N^\prime\psi_N,\psi_N)+
(\B{B}_N^\prime\psi_N,\delta_N)\\
&=(\B{B}_N^\prime\psi_N,\psi_N)-(\B{B}_N^\prime\psi_N,
(\B{I}+\B{A})^{-1}\B{B}_N^\prime\psi_N).
\end{align}
It follows from (\ref{psiN}) that
\[
||\psi_N||^2+(\psi_N,\B{A}\psi_N)-(\psi_N,\B{B}_N^\prime\psi_N)=(\psi_N,e).
\]
Hence, due to the positivity of the operator $\B{A}$,
\begin{equation}                                                 \label{psiNe}
||\psi_N||^2-||\B{B}_N'||\cdot
||\psi_N||^2<|(\psi_N,e)|.
\end{equation}
In the next section we will show that
\begin{equation}                                                 \label{CN}
|(\psi_N,e)|<KN,
\end{equation}
where the constant $K$ does not depend on $x,y,t$ if
$\xi=x-C(Y)t<\frac{1}{2p_0}\log t^{\frac{N+3}{2}+\varepsilon}$.
According to (\ref{bBN'}) we
have $||\B{B}_N'||\to 0$ as $t\to\infty$.
Therefore taking into account
(\ref{psiNe}) and (\ref{CN}) we get
\begin{equation}                                                 \label{2CN}
||\psi_N||^2<2KN
\end{equation}
for sufficiently large $t$. In turn, it follows from (\ref{dele}),
(\ref{2CN}) and (\ref{A<1}) that
\begin{equation}                                                 \label{dNe}
|(\delta_N,e)|\le\frac{K_N(Y)}{t^{1/2-\varepsilon}},
\end{equation}
where $K_N$ does not depend
on $x,y,t$ if $\xi=x-C(Y)t<\frac{1}{2p_0}\log t^{\frac{N+3}{2}+\varepsilon}$.

\begin{rem}                                                       \label{rem.1}
The same estimate is also valid for the derivative
$\frac{\partial}{\partial x}(\psi_N,e)$.
It can be proved using an analytic continuation
of $(\psi_N,e)$ in some strip $|\Im x|<\beta$.
\end{rem}

Thus according to (\ref{u3}) we need to investigate the solution of
(\ref{psiN}). This equation is an integral equation with degenerate kernel:
\begin{align}                                                 \label{AN3}
A_N(p,q,\mu,\nu,x,Y,t)&=\sum_{i,j=0}^N C_{ij}(\lambda-k_0)^i
(\bar k- \bar k_0)^j\notag\\
&\quad
\times E(p,q)E(\mu,\nu)\chi_{N,Y}(p,q) \chi_{N,Y}(\mu,\nu)g(\mu,\nu)
\end{align}
and right-hand side
\begin{equation}                                                 \label{E3}
E(p,q)=E(p,q,x,Y,t)=e^{p(x-f(p,q,Y)t)}.
\end{equation}
Due to the specific form of the kernel we look for a solution of
(\ref{psiN}) in the form:
\begin{equation} 
\label{psiN3}
\psi_N(p,q,x,Y,t)=\sum_{j=0}^N \psi^{(N)}_j(x,Y,t)
(\bar k-\bar k_0)^jE(p,q)\chi_{N,Y}(p,q).
\end{equation}
Substituting (\ref{psiN3}) into (\ref{psiN}) and taking into account 
(\ref{AN3}),
(\ref{E3}) we obtain a system of linear algebraic equations for
the functions
$\psi^{(N)}_j=\psi^{(N)}_j(x,Y,t)$:
\begin{equation} 
	\label{laeq}
\psi^{(N)}_j+\sum_{l=0}^N A_{jl}^{(N)}\psi^{(N)}_l =\delta_{j0}
\qquad j=0,1,\dots,N,
\end{equation}
where $\delta_{00}=1$ and $\delta_{j0}=0$ for $j=1,2,\dots,N$,
\begin{equation} 
\label{AijN}
A^{(N)}_{ij}=A^{(N)}_{ij}(x,Y,t)=
\sum_{l=0}^N C_{il}J_{lj}(x,y,t),
\end{equation}
and the integrals $J_{lj}$ are defined by
\begin{equation}                                        
	\label{Jlj}
J_{lj}(x,Y,t)=\int_{G_{N,Y}} E^2(p,q)(k-k_0)^l
(\bar k-\bar k_0)^jg(p,q)\,d p\,d q.
\end{equation}
The solution of the system (\ref{laeq}) is given by
\begin{equation} 
	\label{det}
\psi^{(N)}_j(x,Y,t)= \frac{D^{(N)}_j(x,Y,t)}{D^{(N)}(x,Y,t)},
\end{equation}
where $D^{(N)}(x,Y,t)=\det[I^{(N)}+A^{(N)}(x,Y,t)]$ is the 
determinant of the matrix with
entries
$\delta_{ij}+A^{(N)}_{ij}(x,Y,t)$ ($i,j=0,1,\dots,N)$, and 
$D_j^{(N)}(x,Y,t)$ is the
determinant of the matrix obtained by replacing the $j$-th column of
$I^{(N)}+A^{(N)}_{ij}$ by the column $(1,0,\dots,0)^{\top}$.

It follows from (\ref{psiN3}) and (\ref{det}) that
\begin{equation}                                                  \label{F/D}
(\psi_N,e)=\frac{F^{(N)}(x,Y,t)}{D^{(N)}(x,Y,t)},
\end{equation}
where $F^{(N)}(x,Y,t)$ is the determinant of the matrix obtained by
replacing the first line
of $I^{(N)}+A^{(N)}(x,Y,t)$ by the line
$\left(J_0(x,Y,t),J_1(x,Y,t),\dots,J_N(x,Y,t)\right)$. Let us now 
note that according to
(\ref{AijN}), $A^{(N)}(x,Y,t)$ is the product of two 
$(N+1)\times(N+1)$ matrices:
\[
A^{(N)}(x,Y,t)=C^{(N)}J^{(N)}(x,Y,t).
\]
$C^{(N)}$ and $J^{(N)}(x,Y,t)$ are the $(N+1)\times(N+1)$ matrices 
with entries $C_{ij}$
and $J_{lj}(x,Y,t)$, respectively. Taking this into account and 
setting $C_{00}$
as a varying
parameter in the matrix $C^{(N)}$ we obtain the determinant formula
\begin{equation} 
\label{Ndet}
(\psi_N,e)=\frac{\partial}{\partial C_{00}}
\log \det[I^{(N)}+A^{(N)}(x,Y,t)].
\end{equation}

\section{Asymptotic behavior of the solution for large time}

First of all let us study the asymptotic behavior of the integrals (\ref{Jlj}).

\begin{lem}                                                      \label{lem.3}
The integrals $J_{ij}(x,Y,t)$ have the following asymptotic
representation:
\[
J_{ij}(x,Y,t)= \frac{g_0|w_0|}{\alpha_0}\frac{h_0^i\bar h_0^j
\Gamma\left(\frac{i+j+3}{2}\right)}{i+j+1}
\left[1+(-1)^{i+j}\right]\frac{e^{2p_0\xi}}{t^{\frac{i+j+3}
{2}}}+\frac{I_{ij}(\xi,Y,t)}{t^{\frac{i+j+4}{2}}}\ .
\]
Here,
\begin{align*}
&g_0=g(k_0),\quad w_0=w(k_0),\\
&h_0=\frac{1}{\alpha_0}\frac{\partial k}{\partial s}(k_0)
= \frac{1}{\alpha_0}\left(\frac{\partial p}{\partial 
s}+i\,\frac{\partial q}{\partial
s}\right)(k_0),
\end{align*}
where $\alpha_0$ is defined in $(\ref{curvat})$,
and the functions $I_{ij}(\xi,Y,t)$ satisfy
\[
|I_{ij}(\xi,Y,t)|\le K_{ij}(Y)(1+|\xi|)e^{2p_0\xi+2\delta_0(N)|\xi|}
\]
with $\delta_0(N)$ defined by $(\ref{delta})$.
\end{lem}

\begin{proof}
Using (\ref{exp}) and taking into account that for $k\in G_{N,Y}$
\[
k-k_0= \frac{\partial k}{\partial s}(k_0)s +
\frac{\partial k}{\partial r}(k_0)r +\mathrm{O}(r^2+|s|^2)
\]
we write down the integral (\ref{Jlj}) in the form:
\begin{align}                                                  \label{Jij4}
\notag
&J_{ij}(x,Y,t)
=g_0|w_0|k_1^i\bar k_1^je^{2p_0\xi}\int_0^{\varepsilon_0(N)}e^{-rt}
\int_{\hat s_-(r)}^{\hat s_+(r)}s^{i+j}[1+ \mathrm{O}(r|s|^\beta)]\,d s\,d r\\
\notag
&\quad
+2g_0|w_0|k_1^i\bar k_1^je^{2p_0\xi}
\int_0^{\varepsilon_0(N)}e^{-rt}
\int_{\hat s_-(r)}^{\hat s_+(r)}[p_1s^{i+j+1}+ p_2rs^{i+j}]
[1+\mathrm{O}(r|s|^\beta)]\,d s\,d r\\
&=J_{ij}^1(x,Y,t)+J_{ij}^2(x,Y,t),
\end{align}
where $\ds k_1=\frac{\partial k}{\partial s}(k_0)$,
$\ds p_1=\frac{\partial p}{\partial s}(k_0)$,
$\ds p_2=\frac{\partial p}{\partial r}(k_0)$
and $\beta=0$ if $i+j=0$ and $\beta=-1$ if $i+j\ge1$.
Integration over $s$ of the first
summand in (\ref{Jij4}) gives the following asymptotic equality:
\[
J^1_{ij}(x,Y,t)=\frac{g_0|w_0|k_1^i\bar k_1^je^{2p_0\xi}} 
{(i+j+1)\alpha_0^{i+j+1}}
\int_0^{\varepsilon_0(N)}
e^{-rt}r^{\frac{i+j+1}{2}}[1+(-1)^{i+j+1}][1+\mathrm{O}(\sqrt r)]\,d r.
\]
Using the asymptotic relation (\ref{gam}) we find
\begin{align}
J_{ij}^1(x,Y,t)
&=\frac{g_0|w_0|}{(i+j+1)\alpha_0} \left(\frac{k_1}{\alpha_0}\right)^i
\left(\frac{\bar k_1}{\alpha_0}\right)^j[1+(-1)^{i+j+1}]\,
\Gamma\left(\frac{i+j+1}{2}\right)\frac{e^{2p_0\xi}} 
{t^{\frac{i+j+3}{2}}}\notag\\
&\quad
+
\mathrm{O}\left(\frac{e^{2p_0\xi}} {t^{\frac{i+j+4}{2}}}\right).
\end{align}
In the same way, taking into account inequality (\ref{E0}), we obtain 
the following
estimate for the second summand in (\ref{Jij4}):
\begin{align}	                                                  \label{J2}
|J_{ij}^2(x,Y,t)|
&\le
K_{ij}'(Y)|\xi|\,e^{2p_0\xi+2\delta_0(N)|\xi|}\int_0^\infty
e^{-rt}r^{\frac{i+j+2}{2}}[1+\mathrm{O}(\sqrt r)]\,d r\notag\\
&\le
K_{ij}(Y)\frac{|\xi|\,e^{2p_0\xi+2\delta_0(N)|\xi|}} {t^{\frac{i+j+4}{2}}},
\end{align}
where the functions $K_{ij}(Y)$ do not depend on $\xi$ and $t$. The 
statement of the lemma
follows from (\ref{Jij4})-(\ref{J2}).
\end{proof}

Now, using (\ref{Ndet}), let us study the large time asymptotic
behavior of the function $D^{(N)}(x,Y,t)=\det [I^{(N)}+A^{(N)}(x,Y,t)]$.

\begin{lem}                                                  \label{lem.4}
We have the following asymptotic relation
\[
\det[I^{(N)}+A^{(N)}(x,Y,t)]=
1+\sum_{n=1}^N\det C^{(n)}\det\Gamma^{(n)}(Y)
\frac{e^{2p_0\xi}}{t^{\frac{n(n+2)}{2}}}[1+\delta_n(\xi,Y,t)]
\]
where $C^{(n)}$ and $\Gamma^{(n)}(Y)$ are $n\times n$ matrices with entries
($i,j=0,1,\dots,n-1$)
\[
{\frac{(i+j)!}{i!j!(2p_0)^{i+j+1}}}\ \mbox{ and }\
{\frac{g_0(Y)|w_0(Y)|}{\alpha_0(Y)}\,\frac{h_0^i(Y) \bar h_0^j(Y)\Gamma\left(
\frac{i+j+3}{2}\right)}{i+j+1}},
\]
respectively, and the functions
$\delta_n(\xi,Y,t)$ satisfy
\begin{equation} 
\label{deln4}
|\delta_n(\xi,Y,t)|\le\frac{K_n(Y)}{t^{1/4}}\quad\text{if }
\xi<\frac{1}{2p_0}\log t^{\frac{N+3}{2}+\varepsilon}\text{ with 
}0<\varepsilon<1/4.
\end{equation}
\end{lem}

\begin{proof}
Let us denote by $\tilde D^{(N)}(x,Y,t;\lambda_0,\dots,\lambda_N)$ 
the determinant of
the matrix $\Lambda^{(N)}+A^{(N)}(x,Y,t)$, where $\Lambda^{(N)}={\rm
diag}(\lambda_0,\dots,\lambda_N)$ is the diagonal matrix depending on 
$N+1$ parameters
$\lambda_0,\dots,\lambda_N$. Clearly,
\[
\tilde D^{(N)}(x,Y,t;1,\dots,1)=
D^{(N)}(x,Y,t)=\det
\left [ I^{(N)}+A^{(N)}(x,Y,t)\right].
\]
This determinant is a polynomial with respect to the $\lambda_k$'s:
\begin{align}                                                 \label{PlN}
&\tilde D^{(N)}(x,Y,t;\lambda_0,\dots,\lambda_N)
=
\lambda_0\dots \lambda_N+
\hat\lambda_0\lambda_1\dots \lambda_ND^{(1)}_0(x,Y,t)\notag\\
&\qquad
+\lambda_0\hat\lambda_1\lambda_2\dots \lambda_ND^{(1)}_1(x,Y,t)+
\dots +\lambda_0\dots \lambda_{N-1}\hat\lambda_ND^{(1)}_N(x,Y,t)\notag\\
&\qquad
+\hat\lambda_0\hat\lambda_1\lambda_2\dots \lambda_ND^{(2)}_{01}(x,Y,t)+
\dots+\lambda_0\dots \hat\lambda_{i_1}\dots \hat\lambda_{i_n}\dots\lambda_N
D^{(n)}_{i_1\dots i_n}(x,Y,t)+\dots\notag\\
&\qquad
+D^{(N)}_{0\dots N}(x,Y,t),
\end{align}
where $D^{(n)}_{i_1\dots i_n}(x,Y,t)$ is the determinant of the 
$n\times n$ matrix
with entries $A_{i_ri_p}(x,Y,t)$ ($r,p=1,\dots,n$);
the hat means that the corresponding summand is absent.
Taking into account
(\ref{AijN}) and using Lemma 4.1 we obtain
\begin{align*}
D^{(n)}_{i_1\dots i_n}(x,Y,t)
&=
\det C^{(n)}\det J^{(n)}(x,Y,t) +d^{(n)}(x,Y,t)\\
&= \det C^{(n)}\det
\Gamma^{(n)}(Y)\frac{e^{2p_0n\xi}} {t^{\frac{n(n+2)}{2}}}+d_1^{(n)}(x,Y,t),
\end{align*}
where
$J^{(n)}(x,Y,t)$ is the $n\times n$ matrix with entries
$J^{(n)}_{ij}(x,Y,t)$ defined by (\ref{Jlj}), $C^{(n)}$ and 
$\Gamma^{(n)}(Y)$ are
defined in Lemma 4.2 ($\det C^{(n)}>0$,  $\det\Gamma^{(n)}(Y)>0$ as
Gram determinants), and the functions $d^{(n)}(x,Y,t)$,
$d_1^{(n)}(x,Y,t)$ satisfy
\begin{equation}                                                 \label{dd1}
|d^{(n)}(x,Y,t)|, \ |d_1^{(n)}(x,Y,t)|<
K_d(Y)\frac{e^{2p_0n\xi}}{t^{\frac{n(n+2)}{2}+1/2}}
\left(1+\frac{e^{2\delta_0(N)|\xi|}}{t^{1/2}}\right).
\end{equation}
The determinants $D^{(n)}_{i_1\dots i_n}$ with 
$i_1+i_2+\dots+i_n>\frac{n(n-1)}{2}$ also
satisfy (\ref{dd1}). Taking all this into account and setting $\lambda_i=1$
($i=0,\dots,N$) in (\ref{PlN}), we obtain the assertion of Lemma \ref{lem.4}.
\end{proof}

\begin{rem}                                                       \label{rem.2}
A more precise analysis of the determinants
$D^{(n)}_{i_1\dots i_n}$ shows that the
derivatives of the functions with respect to $\xi$ and to $C_{00}$ have the
same estimates as in (\ref{deln4}).
\end{rem}

Let us use the following equality, which is proved in \cite{FS}:
\[
n\det C_0^{(n)}=\det C_1^{(n-1)},
\]
where $C_0^{(n)}$ is the $n\times n$ matrix with entries
$\ds\frac{(i+j)!}{i!j!}$ $(i,j=0,\dots,n-1)$, and $C_1^{(n-1)}$ is the
$(n-1)\times(n-1)$ matrix
with entries $\ds\frac{(i+j)!}{i!j!}$ ($i,j=1,\dots,n-1$).
This allows us to obtain the relation:
\begin{equation} 
\label{dC00}
\frac{\partial}{\partial C_{00}}\det C^{(n)}=2p_0n
\det C^{(n)}\!\!\bigm\vert_{\,C_{00}
=(2p_0)^{-1}},
\end{equation}
where $C^{(n)}$ is as above.

Now taking into account (\ref{u3}), (\ref{dNe}), (\ref{Ndet}), Lemma 
\ref{lem.4},
Remarks \ref{rem.1} and \ref{rem.2}
and (\ref{dC00}) we obtain the following asymptotic formula for the solution:
\begin{equation} 
	\label{u4}
u(x,y,t)=2\frac{\partial^2}{\partial\xi^2}\log \left[ 1+\sum_{n=1}^N
\det C^{(n)}\det\Gamma^{(n)}(Y)\frac{e^{2p_0n\xi}}{t^{\frac{n(n+2)}{2}}}
\right]_{\xi=x-C(Y)t}+\mathrm{O}\left(t^{-1/4}\right)
\end{equation}
in the domain
$\{(x,y)\in\D{R}^2\mid x<C(Y)t+\frac{1}{2p_0}\log 
t^{\frac{N+3}{2}+\varepsilon}\}$ as
$t\to\infty$.

\begin{rem}
This asymptotic formula is uniform with respect to $y$ because, in (\ref{dNe})
and (\ref{deln4}), $K_N(Y)$ is uniformly bounded with respect to $y$. 
That follows from the
compactness of the contour $\Gamma$ and the positivity of its curvature.
\end{rem}

To push further the asymptotic analysis of the determinant formula (\ref{u4})
let us introduce notations:
\begin{align}                                                     \label{D_N}
&\Delta_N(\xi,Y,t)=1+\sum_{n=1}^N R_n(Y)\frac{e^{2p_0n\xi}}
{t^{n(n+2)/2}}\\                                                   \label{Rn}
&R_n(Y)=\det C^{(n)}\det \Gamma^{(n)}(Y)=
\frac{\left(\frac{g_0|w_0|}{\alpha_0}\right)^n|h_0|^{n(n-1)}}
{\left(2p_0\right)^{2n+2}\prod\limits_{i=0}^{n-1}(i!)^2} 
\Delta^{(n)}_1\Delta^{(n)}_2,
\end{align}
where $\Delta_1^{(n)}>0$, $\Delta_2^{(n)}>0$ are the determinants of the
$n\times n$ matrices
$\Gamma^{(n)}_1$ and $\Gamma^{(n)}_2$ with entries $\Gamma(i+j+1)$ and
$\Gamma\bigl((i+j+3)/2\bigr)\frac{1+(-1)^{(i+j)}}{i+j+1}$
($i,j=0,1,\dots,n-1$), respectively. They are positive as Gram
determinants.

 From (\ref{u4}) and (\ref{D_N}) it follows that
\begin{equation} 
\label{sim}
u(x,y,t)\sim u_N(x,Y,t)=2\frac{\partial^2}{\partial\xi^2}
\log\Delta_N(\xi,Y,t)\vert_{\xi=x-C(Y)t}=
\frac{\Delta_N^{\prime\prime}\Delta_N-(\Delta_N^\prime)^2} {\Delta_N^2}
\end{equation}
and
\begin{equation} 
\label{sim1}
\Delta_N^{\prime\prime}\Delta_N-(\Delta_N^\prime)^2= 4p_0^2
\sum_{n,l=0}^N\frac{(n-l)^2R_nR_le^{2(n+l)p_0\xi}} 
{t^{\frac{n(n+2)+l(l+2)}{2}}}\ .
\end{equation}
Let us cover the domain $\xi<\frac{1}{2p_0} \log 
t^{\frac{N+3}{2}+\varepsilon}$ by the
intervals
\begin{align*}
&I_1(t)=\left\{ -\infty<\xi<\frac{1}{2p_0}\log t^{2+\varepsilon}\right\},\\
&I_2(t)=\left\{\frac{1}{2p_0}\log t^{2-\varepsilon}<\xi
<\frac{1}{2p_0} \log t^{3+\varepsilon}\right\}\\
&\dots\dots\dots\dots\dots\dots\dots\dots\dots\dots\dots\dots\dots\dots\\
&I_n(t)=\left\{\frac{1}{2p_0}\log t^{n-\varepsilon}<\xi
< \frac{1}{2p_0}\log t^{(n+1)+\varepsilon}\right\},\\
&\dots\dots\dots\dots\dots\dots\dots\dots\dots\dots\dots\dots
\dots\dots\dots\dots\dots\dots\\
&I_{\left [\frac{N+1}{2}\right]}(t)=\left\{\frac{1}{2p_0}\log
t^{\left[\frac{N+1}{2}\right] -\varepsilon}<\xi <\frac{1}{2p_0}
\log t^{\frac{N+3}{2}+\varepsilon} \right\}.
\end{align*}
Taking into account (\ref{sim}) and (\ref{sim1}), we obtain
\[
\Delta_N^2=\left [\frac{R_{n-1}(Y)e^{2(n-1)p_0\xi}} {t^{\frac{(n-1)(n+1)}{2}}}
+ \frac{R_{n}(Y)e^{2np_0\xi}}{t^{\frac{n(n+2)}{2}}}\right ]
\bigl(1+\mathrm{O}\bigl(t^{-1/2}\bigr)\bigr)
\]
and
\[
\Delta_N^{\prime\prime}\Delta_N-(\Delta_N^\prime)^2= 4p_0^2
\frac{2R_n(Y)R_{n-1}(Y)e^{2(2n-1)p_0\xi}} {t^{n(n+2)+(n-1)(n+1)}}\,
\bigl(1+\mathrm{O}\bigl(t^{-1/2}\bigr)\bigr),
\]
as $\xi\in I_n(t)$ and $t\to\infty$. Hence, by virtue of (\ref{D_N}) 
and (\ref{sim})
we obtain
\begin{align}	                                                 \label{ufin}
u(x,y,t)
&=\sum_{n=1}^{\left[\frac{N+1}{2}\right]}
\frac {2p_0^2(Y)}{\cosh^2\left [p_0(Y)
\left (x-C(Y)t+ \frac{1}{2p_0(Y)}\log t^{n+1/2}+x_n^0(Y)\right 
)\right]}\notag\\
&\quad
+\mathrm{O}\left (t^{-1/4}\right),
\end{align}
where $x^0_n(Y)=\frac{1}{2p_0(Y)}\log\frac{R_n(Y)}{R_{n-1}(Y)}$ and $Y=y/t$.

\begin{rem}                                                    \label{rem.4}
The estimate (\ref{CN}) follows from (\ref{u2}),
(\ref{u4}) and (\ref{ufin}).
\end{rem}

Thus we have proved the following result:

\begin{thm}
Let the contour $\Gamma$ be compact, of class $C^2$ on $\bar\Omega$,
and with everywhere
positive curvature.
Assume the function $f(p,q,Y)$ attains its minimal value at a unique
point
$k_0(Y)=p_0(Y)+i q_0(Y)\in\Gamma$, for any $Y=y/t$.

Then the solution
$u(x,y,t)$ of the
KP-I equation defined everywhere by $(\ref{psieq})$ and $(\ref{u2})$ in
the domain
\[
D_N=\Bigl\{(x,y)\mid-\infty<y<\infty,\,x<C(Y)t+\frac{1}{2p_0(Y)}
\log t^{[(N+1)/2]+1+\varepsilon}\Bigr\}\quad (0<\varepsilon<1/4)
\]
has the asymptotic behavior
defined by $(\ref{u4})$ and $(\ref{ufin})$ for $t\to\infty$.
\end{thm}


\label{lastpage}
\end{document}